\documentclass[12pt,3p,times]{elsarticle}

\usepackage{enumitem}
\usepackage{hyperref}
\usepackage{amsmath}
\usepackage{pifont}
\usepackage[utf8]{inputenc}
\makeatletter
\def\LaTeX{\leavevmode L\raise.42ex
    \hbox{\kern-.3em\size{\sf@size}{0pt}\selectfont A}\kern-.15em\TeX}
\makeatother

\newcommand{\BibTeX}{{\rm B\kern-.05em{\sc
          i\kern-.025emb}\kern-.08em\TeX}}

\makeatletter
\def\@currentlabel{2.1}\label{e:dispaa}
\def\@currentlabel{2.21}\label{e:dispau}
\def\@currentlabel{2.22}\label{e:dispav}
\def\@currentlabel{2.23}\label{e:dispaw}
\def\@currentlabel{2.24}\label{e:dispax}
\def\theequation{\thesection.\@arabic\c@equation}
\makeatother

\renewcommand{\theequation}{\arabic{section}.\arabic{equation}}

\newcommand{\R}{\mathbb R}

\def \O{\Omega}
\everymath{\displaystyle}

\newtheorem{thm}{Theorem} [section]
\newtheorem{lem}{Lemma} [section]
\newtheorem{prop}{Proposition} [section]

\newtheorem{definition}{Definition} [section]

\newtheorem{rem}{Remark}[section]

\renewcommand{\theequation}{\thesection.\arabic{equation}}
\renewcommand{\thesection}{\arabic{section}}
\renewcommand{\theequation}{\thesection.\arabic{equation}}
\let\ssection=\section\renewcommand{\section}{\setcounter{equation}{0}\ssection}
\begin{document}
\begin{frontmatter}
\title{{\bf{New class of sixth-order nonhomogeneous $p(x)$-Kirchhoff problems with sign-changing
weight functions}}}
\author[mk0,mk1,mk2]{Mohamed Karim Hamdani}
\ead{hamdanikarim42@gmail.com}
\cortext[cor1]{Corresponding author\; Du\u{s}an D. Repov\v{s}: dusan.repovs@guest.arnes.si}
\author[dr0]{Nguyen Thanh Chung}
\ead{ntchung82@yahoo.com}
\author[dr1,dr2]{ Du\u{s}an D. Repov\v{s}\corref{cor1}}
\ead{dusan.repovs@guest.arnes.si}
\begin{center}
\address[mk0] {Science and Technology for Defense Laboratory  LR19DN01, Military Research Center, Aouina, Tunisia.}
\address[mk1]{Military School of Aeronautical Specialities, Sfax, Tunisia.}
\address[mk2]{Mathematics Department, University of Sfax, Faculty of Science of Sfax, Sfax, Tunisia.}
\address[dr0]{Department of Mathematics, Quang Binh University, 312 Ly Thuong Kiet, Dong Hoi, Quang Binh, Vietnam.}
\address[dr1] {Faculty of Education and Faculty of Mathematics and Physics, University of Ljubljana, Ljubljana, Slovenia}
\address[dr2] {Institute of Mathematics, Physics and Mechanics,  Ljubljana, Slovenia}
\end{center}
\begin{abstract}In this paper, we prove the existence of multiple solutions for the following sixth-order $p(x)$-Kirchhoff-type
problem
\begin{eqnarray*} \label{abs10}
\begin{cases}
-M\left( \int_\O \frac{1}{p(x)}|\nabla \Delta u|^{p(x)}dx\right)\Delta^3_{p(x)} u = \lambda f(x)|u|^{q(x)-2}u + g(x)|u|^{r(x)-2}u + h(x) &\mbox{in}\quad \Omega, \\
u = \Delta u = \Delta^2 u = 0, \quad &\mbox{on}\quad \partial\Omega,
\end{cases}
\end{eqnarray*}
where $\O\subset \R^N$ is a smooth bounded domain, $N > 3$, $\Delta_{p(x)}^3u := \operatorname{div}\Big(\Delta(|\nabla \Delta u|^{p(x)-2}\nabla \Delta u)\Big)$ is the $p(x)$-triharmonic operator, $p,q,r \in C(\overline\O)$, $1< p(x) < \frac N3$ for all $x\in \overline\O$, $M(s) = a - bs^\gamma$, $a,b,\gamma>0$, $\lambda>0$, $g: \O \times \R \to \R$ is a nonnegative continuous function while $f,h : \O \times \R \to \R$ are sign-changing continuous functions in $\O$. To the best of our knowledge, this paper is one of the first contributions to the study of the sixth-order $p(x)$-Kirchhoff type problems with sign changing Kirchhoff functions.
\end{abstract}
\end{frontmatter}
\section{\bf Introduction}
Let $\O\subset \R^N$ be a smooth bounded domain and $N > 3$. This paper deals mainly with the following sixth-order $p(x)$-Kirchhoff-type problem
\begin{eqnarray} \label{high}
\begin{cases}
-M\left( \int_\O \frac{1}{p(x)}|\nabla \Delta u|^{p(x)}dx\right)\Delta^3_{p(x)} u = \lambda f(x)|u|^{q(x)-2}u + g(x)|u|^{r(x)-2}u + h(x) &\mbox{in}\quad \Omega, \\
u = \Delta u = \Delta^2 u = 0, \quad &\mbox{on}\quad \partial\Omega,
\end{cases}
\end{eqnarray}
where $p,q,r \in C(\overline\O)$, $1< p(x) < \frac N3$ for all $x\in \overline\O$, $M(s) = a - bs^\gamma$, $a,b,\gamma>0$, $\lambda>0$, $g: \O \times \R \to \R$ is a nonnegative continuous function,  $f,h : \O \times \R \to \R$ are assumed to be continuous functions which may change sign in $\O$, and $$\Delta_{p(x)}^3u := \operatorname{div}\Big(\Delta(|\nabla \Delta u|^{p(x)-2}\nabla \Delta u)\Big)$$ is the $p(x)$-triharmonic operator which is not homogeneous and is related to the variable exponent Lebesgue space $L^{p(x)}(\O)$ and the variable exponent Sobolev space $W^{1,p(x)}(\O)$. These facts imply some difficulties. For example, some classical theories and methods, including the Lagrange multiplier theorem and the theory of Sobolev spaces cannot be applied.

Such problems are called nonlocal problems because of the presence of the function $M$, which implies that the equation contains an integral over $\O$, and is no longer pointwise identity. This causes some mathematical difficulties which make the study of such a problem particularly interesting. We call \eqref{high} a sixth order Kirchhoff type equation because it is related to the stationary analog of the equation
\begin{equation} \label{kirchhoff}
\rho\frac{\partial ^2u}{\partial t^2}-\Bigg(\frac{p_0}{h}+\frac{E}{2L} \int_0^L\left|\frac{\partial u}{\partial x}\right|^2\,dx\Bigg) \frac{\partial ^2u}{\partial x^2}=0,
\end{equation}
where $\rho$, $p_0$, $h$, $E$, $L$ are constants which represent some physical meanings respectively. Eq.~\eqref{kirchhoff} extends
the classical D'Alembert wave equation by considering the effects of the changes in the length of the strings during the vibrations.

This kind of nonlocal problem also appears in other fields, for example in nonlinear elasticity theory and in modelling electrorheological fluids \cite{RR2, R3} and from the study of electromagnetism and elastic mechanics \cite{HR,Z}, and raises many difficult mathematical problems. After this pioneering models, many other applications of differential operators with variable exponents have appeared in a large range of fields, such as image restoration and image processing \cite{CLR,LLP}. We refer the reader to \cite{AM1,H,MR06} for an overview of references on this subject.

Throughout this paper, unless otherwise stated, we shall always assume that exponent $p(x)$ is continuous on $\overline{\O}$ with
$$p_-:=\inf_\O p(x) \leq p(x)\leq p_+:=\sup_\O p(x)<\frac N3$$
and $p^*(x)$ denotes the critical variable exponent related to $p(x)$, defined for all $x\in\overline{\O}$ by the pointwise relation $p^*(x) = \frac{Np(x)}{N-3p(x)}$. In the following, we denote by $[W_0^{3,{p(\cdot)}}(\O)]'$ the dual space of $W_0^{3,{p(\cdot)}}(\O)$ and $q^*(x) = \frac{p^*(x)}{p^*(x)-1}$ the conjugate exponent of $p^*(x)$.

Now, we introduce some conditions for problem \eqref{high} as follows:
\begin{itemize}
\item[$(H_1):$] $1< q(x) < p_- \leq p_+ < (\gamma+1)p_- \leq (\gamma+1) p_+ < r(x) < p^*(x)$ for all $x \in \overline\O$, where $\gamma$ is a positive constant;
\item[$(H_2):$] $f\in L^{q_0(\cdot)+\eta}(\O)$, $0 \leq g \in L^{r_0(\cdot)+\mu}(\O)$ and $h \in L^{\frac{p^*(\cdot)}{p^*(\cdot)-1}} (\O) \cap L^\infty(\O)$ for all $x \in \overline\O$, with $$q_0(x) = \frac{p^*(x)}{p^*(x)-q(x)}, \ r_0(x) = \frac{p^*(x)}{p^*(x)-r(x)}$$ where $\eta,\mu$ are small positive numbers.
\item[$(H'_2):$] $f\in L^{q_0(\cdot)}(\O)\cap L^{\infty}(\O)$ and $g \in L^{r_0(\cdot)}(\O)\cap L^{\infty}(\O)$ with $$q_0(x) = \frac{p^*(x)}{p^*(x)-q(x)}, \  r_0(x) = \frac{p^*(x)}{p^*(x)-r(x)}.$$ Furthermore, there
exists non-empty open domain $\O_0\subset \O$ such that $g(x) > 0$ in $\O_0$.
\end{itemize}

In recent years, great attention has been paid to the study of Kirchhoff problems. This brought new mathematical difficulties that made the study of Kirchhoff type equations particularly interesting. A typical prototype for $M$, due to Kirchhoff in $1883$, is given by
\begin{eqnarray}\label{hhj}
M(t)=a+b\tau^{\alpha-1},\quad a,b\geq 0,\;a+b>0,t\geq 0,
\end{eqnarray}
and
\begin{eqnarray*}
\begin{cases}
\alpha\in(1,+\infty) & \mbox{if } b>0 \\
\alpha =1, & \mbox{if } b=0,
\end{cases}
\end{eqnarray*}
when $M(t) > 0$ for all $t\geq 0$, Kirchhoff problems are said to be nondegenerate and this happens for example if $a>0$ and $b\geq0$ in the model case \eqref{hhj}. Otherwise, if $M(0) = 0$ and $M(t) > 0$ for all $t >0$, the Kirchhoff problems are called degenerate and this occurs in the model case \eqref{hhj} when $a = 0$ and $b>0$. In particularly, Chen-Kuo-Wu in \cite{CKW} studied the following semilinear boundary problem and proved by the Nehari manifold and fibering
maps, the existence of multiple positive solutions
\begin{eqnarray*}
\begin{cases}
-\left( a+b\int_\O |\nabla u|^{2}dx\right)\Delta u = \lambda f(x)|u|^{q-2}u + g(x)|u|^{r-2}u  &\mbox{in}\quad \Omega, \\
u = 0, \quad &\mbox{on}\quad \partial\Omega,
\end{cases}
\end{eqnarray*}
where $\O$ is a smooth bounded domain in $\R^N$, with $1 < q < 2 < p < 2^*=\frac{2N}{N-2}$ and the parameters $a, b, \lambda > 0$. The functions $f(x), g(x) \in C(\overline{\O})$ may change sign on $\O$.

The study of Kirchhoff type equations has already been extended to the case involving the $p$-Laplacian operator of the following form
\begin{eqnarray} \label{p-lap}
\begin{cases}
-M\left( \int_\O |\nabla u|^{p}dx\right)\Delta_{p} u = \lambda f(x)|u|^{q-2}u + g(x)|u|^{r-2}u + h(x) &\mbox{in}\quad \Omega, \\
u = 0, \quad &\mbox{on}\quad \partial\Omega.
\end{cases}
\end{eqnarray}
Chen-Huang-Liu \cite{CHL}
studied the nonhomogeneous case of \eqref{p-lap} (that is $h(x)\neq 0$), $\lambda > 0$, $M(s) = a + bs^k$, $a, b > 0$, $k \geq 0$, $f(x), g(x)$ and $h(x)$ are continuous functions which may change sign on $\O$. The parameters $p, q,r$
satisfy $1 < q < p(k + 1) < r < p^*=\frac{Np}{N-p}$. Using the Mountain pass theorem and Ekeland's variational principle, they showed that
problem \eqref{p-lap} has at least two positive solutions when $\lambda$ is small enough.

For $h(x)\equiv 0$, Huang-Chen-Xiu \cite{HCX} studied problem \eqref{p-lap} where $M(s) = a + bs^k$, $1 < q < p < r < p^*$, and proved that the problem has at least one positive solution when $r > p(k + 1)$ and the functions $f(x), g(x)$ are nonnegative. Motivated by \cite{HCX}, Li-Mei-Zhang \cite{LMZ} considered $M(s) = a + bs$, $1< q < p < r \leq p^*$ and they proved the existence of multiple nontrivial nonnegative solutions by using the Nehari manifold when the weight functions $f(x), g(x)$ change their signs (see also \cite{H2}).

However, many papers generalized the constant case to include the $p(x)$-Laplacian operator, e.g., in \cite{Chung2012}, using variational
methods, we investigated a nonlocal $p(x)$-Laplacian Dirichlet problem
\begin{eqnarray} \label{p(x)-lapl}
\begin{cases}
-M\left( \int_\O \frac{1}{p(x)}|\nabla u|^{p(x)}dx\right)\Delta_{p(x)} u = K(x,u) &\mbox{in}\quad \Omega, \\
u = 0 \quad &\mbox{on}\quad \partial\Omega,
\end{cases}
\end{eqnarray}
and we showed via the mountain pass theorem combined with the Ekeland variational principle the existence of at least two distinct, non-trivial weak solutions in the case that
$$K(x,u)=\lambda\Big(a(x)|u|^{\alpha(x)-2}u+b(x)|u|^{\beta(x)-2}u\Big),$$
where $\lambda$ is a parameter and $a(x),b(x),\alpha(x)$ and $\beta(x)$ satisfy suitable hypotheses and under some suitable conditions on $M$. 

Recently, Hamdani et al. in \cite{HHMR} studied \eqref{p(x)-lapl} when $M(s)=a-bs$ and $K(x,u)=\lambda |u| ^{p(x)-2}u+g(x,u)$, where $\lambda$ is a real parameter, $a, b > 0$ are constants and $g$ is a continuous function satisfies the classical $(AR)$ condition. If $1<p^-<p(x)<p^+<2p^-<q^-<q(x)<p^*(x)$ then the authors proved the existence and multiplicity results via the Mountain Pass theorem and the Fountain theorem. Also for further studies on this subject, we refer the reader to \cite{JF, Chung13, C182,CN10,  DH, LRDZ, MR10, Y}.

 The problems involving $p(x)-$biharmonic operators have been widely investigated. For example, via  the Mountain pass theorem, El Amrouss et al. \cite{AMM} obtained the existence and multiplicity of solutions for  a class of $p(x)-$biharmonic equation of the form
\begin{eqnarray}\label{biharmonic}
\begin{cases}
\Delta^2_{p(x)}u=\lambda |u|^{p(x)-2}u + f(x,u) & \mbox{ in }\quad \Omega, \\
u=\Delta u=0 & \mbox{ on }\quad \partial\Omega,
\end{cases}
\end{eqnarray}
where $\O$ is a bounded domain in $\R^N$, with smooth boundary $\partial \O$, $\lambda\leq 0$ and $f$ satisfies the $(AR)$ condition. On the other hand, similar variational methods are also used to study the $p(x)$-biharmonic operator. For example, see \cite{AMC14, BR16,BRR16,  C181, MAE, MTTL, Zhou} and \cite{CP,HCA,HZCR} for a general Kirchhoff problems with or without the $(AR)$ condition. However, in literature the only result involving a sixth-order problem like \eqref{high} by using variational methods can be found in \cite{Rahal}. Motivated by the above papers, we consider the problem \eqref{high} with the new general nonlocal term $M(s) = a - bs^\gamma$, $a,b,\gamma>0$ and with sign-changing weight functions which presents interesting difficulties to discus the existence of multiple solutions.

Before stating our results, we give the definition of weak solutions for problem \eqref{high}. For this purpose, we denote by $X$ the space $W_0^{1,p(\cdot)}(\O) \cap W_0^{3,p(\cdot)}(\O)$ and define the norm $\|.\|_X$ of $X$ by the formula
$$
\|u\|_X = \|u\|_{1,p(\cdot)} + \|u\|_{2,p(\cdot)} + \|u\|_{3,p(\cdot)}.
$$

It is well known that if $1<p^- \leq p^+ < \infty$ then $(X, \|.\|_X)$ is a separable and reflexive Banach space. Moreover, $\|u\|_X$ and $|\nabla \Delta u|_{p(\cdot)}$ are two equivalent norms on $X$, see \cite{FS,FZ1,KR}.

\begin{definition}\label{def1.1}
We say that $u \in X$ is a weak solution of problem \eqref{high} if $u$ satisfies
\begin{gather*}
\left(a-b\left(\int_\O \frac{|\nabla \Delta u|^{p(x)}}{p(x)}\,dx\right)^{\gamma}\right)\int_\O|\nabla \Delta u|^{p(x)-2}\nabla\Delta u \nabla \Delta v\,dx - \lambda\int_\O f(x)|u|^{q(x)-2}uv\,dx \\
- \int_\O g(x)|u|^{r(x)-2}uv\,dx - \int_\O h(x)v\,dx = 0,
\end{gather*}
for all $v\in X$.
\end{definition}
The main results of this paper are the following:
\begin{thm}\label{the1.1}
Assume that the conditions $(H_{1})$ and $(H_{2})$ hold. Then there exist $\overline\lambda, \delta>0$ such that, for each $\lambda \in (0,\overline\lambda)$, problem \eqref{high} admits at least two nontrivial weak solutions in $X$ provided that $|h|_{\frac{p^*(\cdot)}{p^*(\cdot)-1}} < \delta$.
\end{thm}
\begin{thm}\label{the1.2}
Let $h(x)\equiv 0$ for all $x\in \O$ and assume that the conditions $(H_{1})$ and $(H'_{2})$ hold. Then there exists $\overline\lambda>0$ such that, for each $\lambda \in (0,\overline\lambda)$, problem \eqref{high} admits at least two nontrivial weak solutions in $X$.
\end{thm}

The paper is organized as follows. In Section \ref{sect1}, we give the notations and recall some useful lemmas concerning the variable exponent Lebesgue and Sobolev spaces. In Section \ref{sect3}, we give some lemmas which are important for the proofs of our main results. In Section \ref{sect4}, we prove Theorem \ref{the1.1} (we omit the proof of Theorem \ref{the1.2} since it is very similar).

\section{\bf Variable exponent Lebesgue and Sobolev spaces}\label{sect1}
For the convenience of the reader, we recall in this section some results concerning spaces $L^{p(\cdot)}(\O)$ and $W^{r,p(\cdot)}(\O)$ which
we call generalized Lebesgue-Sobolev spaces. Denote
$$C_+(\overline{\O})=\left\{p(x):\; p(x)\in C(\overline{\O}),\; p(x)>1, \;\hbox{for  all} \;  x\in \overline{\O}\right\}.$$
For any $p(x)\in C_+(\overline{\O})$, we introduce the variable exponent Lebesgue space
$$L^{p(\cdot)}(\O)=\left\{u(x): u(x) \mbox{ is a measurable real-valued function such that } \int_\O |u(x)|^{p(x)}\,dx<\infty\right\},$$
endowed with the so-called Luxemburg norm
$$\|u\|_{L^{p(\cdot)}(\O)}=|u|_{p(\cdot)}=\inf \left\{\mu>0;\int_\O\left|\frac{u(x)}{\mu}\right|^{p(x)}\,dx\leq 1\right\},$$
which is a separable and reflexive Banach space. A thorough variational analysis of the problems with variable exponents has been developed in the monograph by R\u{a}dulescu and Repov\v{s} \cite{RD} (we refer the reader also to \cite{Diening,KR}).

\begin{prop}[see \cite{Y}]
 The space $(L^{p(\cdot)}(\O), |\cdot|_{p(\cdot)})$ is separable, uniformly convex, reflexive and its conjugate space is $L^{p(\cdot)}(\O), |\cdot|_{q(\cdot)}$ where $q(x)$ is the conjugate function of $p(x)$ i.e
$$\frac{1}{p(x)}+\frac{1}{q(x)}=1,\;\;\hbox{for all} \  x\in \O.$$
For all $u\in L^{p(\cdot)}(\O)$ and $v\in L^{q(\cdot)}(\O)$ the H\"older type inequality
\begin{eqnarray*}
  \left|\int_\O uvdx\right|\leq\left(\frac{1}{p^-}+\frac{1}{q^-}\right)|u|_{p(\cdot)}|v|_{q(\cdot)} \leq 2|u|_{p(\cdot)}|v|_{q(\cdot)} \end{eqnarray*}
  holds.
\end{prop}

The inclusion between Lebesgue spaces also generalizes the classical framework, namely if $0<|\O|<\infty$ and $p_1$, $p_2$ are variable exponents such that $p_1 \leq p_2$ in $\O$, then there exists a continuous embedding $L^{p_2(\cdot)}(\O)\to L^{p_1(\cdot)}(\O)$. An important role in manipulating the generalized Lebesgue-Sobolev spaces is played by the $p(\cdot)-$modular of the $L^{p(\cdot)}(\O)$ space, which is the modular $\rho_{p(\cdot)}$ of the space $L^{p(\cdot)}(\O)$
\begin{equation*}
\rho_{p(\cdot)}(u):=   \int_{\Omega} |u|^{p(x)} \,dx.
\end{equation*}

\begin{lem}\label{lemmaineq}
If $u_n, u \in L^{p(\cdot)}$ and $p_{+} < +\infty$, then the following
properties hold:
\begin{enumerate}
\item $|u|_{p(\cdot)} > 1 \Rightarrow
|u|_{p(\cdot)}^{p_{-}} \leq \rho_{p(\cdot)}(u) \leq |u|_{p(\cdot)}^{p_{+}}$;

\item $|u|_{p(\cdot)} < 1 \Rightarrow   |u|_{p(\cdot)}^{p_{+}}
\leq \rho_{p(\cdot)}(u) \leq |u|_{p(\cdot)}^{p_{-}}$;

\item $|u|_{p(\cdot)} < 1$ (respectively $= 1; > 1) \Longleftrightarrow
 \rho_{p(\cdot)}(u) < 1$ (respectively $= 1; > 1$);\label{gggg}

\item $|u_n|_{p(\cdot)} \to 0$ (respectively
$\to +\infty) \Longleftrightarrow \rho_{p(\cdot)}(u_n) \to 0$
 (respectively $\to +\infty$);
\item $\lim_{n\to \infty}|u_n-u|_{p(\cdot)}=0 \Longleftrightarrow \lim_{n\to \infty}\rho_{p(\cdot)}(u_n-u)=0$.
\end{enumerate}
\end{lem}

The Sobolev space with variable exponent $W^{r,p(\cdot)}(\O)$ is defined as
\begin{equation*}
   W^{r,p(\cdot)}(\Omega):=\Big\{u \in L^{p(\cdot)}(\Omega):
D^\alpha u \in L^{p(\cdot)}(\Omega),\;|\alpha|\leq r \Big\},
\end{equation*}
where $D^\alpha u=\frac{\partial^{|\alpha|}}{\partial x_1^{\alpha_1}\partial x_2^{\alpha_2}...\partial x_N^{\alpha_N}u}$, with $\alpha=(\alpha_1,...\alpha_N)$ is a multi-index and $|\alpha|=\sum_{i=1}^{N}\alpha_i.$ The space $W^{r,p(\cdot)}(\Omega)$ is a reflexive and separable Banach space if $1<p_-\leq p_+<+\infty$ and equipped with the norm
$$\|u\|_{r,p(\cdot)}:=\sum_{|\alpha|\leq r}|D^\alpha u|_{p(\cdot)}.$$

Let $W_0^{r,p(\cdot)}(\O)$ denote the completion of $C_0^\infty (\O)$ in $W^{r,p(\cdot)}(\O)$. As shown in (\cite{Diening},
Corollary $11.2.4$), the space $W_0^{r,p(\cdot)}(\O)$ coincides with the closure in $W^{r,p(\cdot)}(\O)$ of the set of all $W^{r,p(\cdot)}(\O)$-functions with compact support.

\begin{prop}[see \cite{FZ1}]\label{embedding} Assume that $s\in C_+(\overline{\O})$ satisfies $s(x) \leq p^*(x)$ for all $x\in \overline{\O}$. Then there is a continuous embedding $X \hookrightarrow L^{s(\cdot)}(\O)$. If we replace $\leq$ with $<$, then this embedding is compact.
\end{prop}

In the light of the variational structure of \eqref{high}, we look for critical points of the associated Euler-Lagrange functional $J:X\to \R$ defined as
\begin{gather}\label{e3.1}
J(u) = a\int_\O \frac{|\nabla \Delta u|^{p(x)}}{p(x)}\,dx-\frac{b}{\gamma+1}\left(\int_\O \frac{|\nabla \Delta u|^{p(x)}}{p(x)}\,dx\right)^{\gamma+1} - \lambda\int_\O \frac{f(x)}{q(x)}|u|^{q(x)}\,dx \notag\\
- \int_\O \frac{g(x)}{r(x)}|u|^{r(x)}\,dx - \int_\O h(x)u\,dx,\;\hbox{for all} \  u\in X.
\end{gather}

\noindent 
Note that $J$ is a $C^1(X, \R)$ functional and
\begin{eqnarray} \label{deriv-fun}
\langle J'(u), {v}\rangle &= \left[a-b\left(\int_\O \frac{|\nabla \Delta u|^{p(x)}}{p(x)}\,dx\right)^{\gamma}\right]\int_\O|\nabla \Delta u|^{p(x)-2}\nabla\Delta u \nabla \Delta v\,dx - \lambda\int_\O f(x)|u|^{q(x)-2}uv\,dx \notag\\
&\qquad\qquad- \int_\O g(x)|u|^{r(x)-2}uv\,dx - \int_\O h(x)v\,dx = 0
\end{eqnarray}
for any $v\in X$. Thus, critical points of $J$ are weak solutions of \eqref{high}.

\section{\bf  Some Lemmas}\label{sect3}
In order to prove our main result - Theorem \ref{the1.1} - we need to apply the Mountain pass theorem and the Ekeland variational principle. We first prove the following lemmas.

\begin{lem}\label{lem1}
Assume that $f$ satisfies $(H_1)-(H_2)$. Then there exist $\overline\lambda, \delta, \rho, \alpha>0$ such that for $\lambda \in (0, \overline \lambda)$ and $|h|_{\frac{p^*(\cdot)}{p^*(\cdot)-1}}<\delta$, we have $J(u) \geq \alpha$ for all $u\in X$ with  $\|u\|_X = \rho$. Moreover, there exists $e\in X$ with $\|e\|_X > \rho$, such that $J(e)<0$.
\end{lem}

{\bf Proof.} {\bf Step 1.} From $(H_1)$ and $(H_2)$, we can note that $f\in L^{q_0(\cdot)+\eta}(\O)$, $g \in L^{r_0(\cdot)+\mu}(\O)$ imply $f\in L^{q_0(\cdot)}(\O)$, $g \in L^{r_0(\cdot)}(\O)$. Then by Proposition \ref{embedding}, there exist constants $C_1, C_2, C_3>0$ such that
\begin{eqnarray}\label{eer1}
\int_\O |f(x)||u|^{q(x)}\,dx & \leq & 2|f|_{q_0(\cdot)}\left||u|^{q(x)}\right|_{\frac{p^*(x\cdot)}{q(\cdot)}} \notag\\
& \leq & C_1|f|_{q_0(\cdot)}\max\left\{\|u\|^{q_+}_X, \|u\|^{q_-}_X\right\},
\end{eqnarray}
\begin{eqnarray}\label{er2}
\int_\O |g(x)||u|^{r(x)}\,dx & \leq & 2|g|_{r_0(\cdot)}\left||u|^{r(x)}\right|_{\frac{p^*(\cdot)}{r(\cdot)}} \notag\\
& \leq & C_2|g|_{r_0(\cdot)}\max\left\{\|u\|^{r_+}_X, \|u\|^{r_-}_X\right\}
\end{eqnarray}
and by Young's inequality,
\begin{eqnarray}\label{er3}
\int_\O |h(x)||u|\,dx & \leq & 2|h|_{\frac{p^*(\cdot)}{p^*(\cdot)-1}}|u|_{p^*(\cdot)} \notag\\
& \leq & C_3|h|_{\frac{p^*(\cdot)}{p^*(\cdot)-1}} \|u\|_X \notag\\
& \leq & \varepsilon C_3\|u\|^{p_+}_X + C_3C_\varepsilon |h|_{\frac{p^*(\cdot)}{p^*(\cdot)-1}}^{\frac{p_+}{p_+-1}}
\end{eqnarray}
for all $u\in X$ and given $\varepsilon>0$, $C_\varepsilon$ is a positive constant depending on $\varepsilon$.

Now, for any $u\in X$ with $\|u\|_X < 1$, from \eqref{eer1}-\eqref{er3} we get
\begin{eqnarray*}\label{er4}
J(u) & = & a\int_\O \frac{|\nabla \Delta u|^{p(x)}}{p(x)}\,dx-\frac{b}{\gamma+1}\left(\int_\O \frac{|\nabla \Delta u|^{p(x)}}{p(x)}\,dx\right)^{\gamma+1} - \lambda\int_\O \frac{f(x)}{q(x)}|u|^{q(x)}\,dx \notag\\
&  & - \int_\O \frac{g(x)}{r(x)}|u|^{r(x)}\,dx - \int_\O h(x)udx \notag\\
& \geq & \frac{a}{p_+}\int_\O|\nabla \Delta u|^{p(x)}\,dx -\frac{b}{\gamma+1}\left(\int_\O \frac{|\nabla \Delta u|^{p(x)}}{p(x)}\,dx\right)^{\gamma+1} - \frac{\lambda}{q_-}\int_\O |f(x)||u|^{q(x)}\,dx \notag\\
& & - \frac{1}{r_-}\int_\O |g(x)||u|^{r(x)}\,dx - \int_\O |h(x)||u|\,dx \notag\\
& \geq & \frac{a}{p_+}\|u\|^{p_+}_X - \frac{b}{(p_-)^{\gamma+1}(\gamma+1)}\|u\|^{p_-(\gamma+1)}_X-\frac{\lambda C_1}{q_-}|f|_{q_0(\cdot)}\|u\|^{q_-}_X-\frac{C_2}{r_-}|g|_{r_0(\cdot)}\|u\|^{r_-}_X \notag\\
& & - \varepsilon C_3\|u\|^{p_+}_X - C_3C_\varepsilon |h|_{\frac{p^*(\cdot)}{p^*(\cdot)-1}}^{\frac{p_+}{p_+-1}}.
\end{eqnarray*}
Choosing $\varepsilon = \frac{a}{2p_+C_3}$, this leads to
\begin{eqnarray*}\label{er5}
J(u) &\geq & \frac{a}{2p_+}\|u\|^{p_+}_X - \frac{b}{(p_-)^{\gamma+1}(\gamma+1)}\|u\|^{p_-(\gamma+1)}_X - C_1|f|_{q_0(\cdot)}\|u\|^{q_-}_X-C_2|g|_{r_0(\cdot)}\|u\|^{r_-}_X \notag\\
& & - C_3C_\varepsilon |h|_{\frac{p^*(\cdot)}{p^*(\cdot)-1}}^{\frac{p_+}{p_+-1}} \notag\\
& = & \|u\|^{p_+}_X\left(\frac{a}{2p_+}-\frac{b}{(p_-)^{\gamma+1}(\gamma+1)}\|u\|^{p_-(\gamma+1)-p_+}_X-\frac{\lambda C_1}{q_-}|f|_{q_0(\cdot)}\|u\|^{q_--p_+}_X-\frac{C_2}{r_-}|g|_{r_0(\cdot)}\|u\|^{r_--p_+}_X\right) \notag\\
& & - C_3C_\varepsilon |h|_{\frac{p^*(\cdot)}{p^*(\cdot)-1}}^{\frac{p_+}{p_+-1}}.
\end{eqnarray*}

Since $p^+ < p_-(\gamma+1) \leq p_+(\gamma+1) < r_-$, we can choose $\rho>0$ sufficiently small so that the following holds
$$
C_\rho = \frac{a}{2p_+} - \frac{b}{(p_-)^{\gamma+1}(\gamma+1)}\rho^{p_-(\gamma+1)-p_+} - \frac{C_2}{r_-}|g|_{r_0(\cdot)}\rho^{r_--p_+} > 0.
$$
Hence, let us choose $\overline \lambda = \frac{C_\rho q_-}{2C_1|f|_{q_0(\cdot)}\rho^{q_--p_+}} > 0$ and $\delta = \frac{1}{2} \left(\frac{C_\rho \rho^{p_+}}{2C_3C_\varepsilon}\right)^{\frac{p_+}{p_+-1}} > 0$. It follows that for each $\lambda \in (0, \overline\lambda)$ and $|h|_{\frac{p^*(\cdot)}{p^*(\cdot)-1}} < \delta$, we have $J(u) \geq \frac{1}{4}C_\rho\rho^{p_+} = \alpha>0$.
\\\\
{\bf Step 2.} Let $\phi_0 \in C_0^\infty(\O_0)$, where $\O_0 \subset \left\{x\in \Omega:~ g(x)>0\right\}$. According to the conditions $(H_1)$ and $(H_2)$, for $t>1$ large enough we have
\begin{eqnarray*}\label{er4}
J(t\phi_0) & = & a\int_\O \frac{|\nabla \Delta t\phi_0|^{p(x)}}{p(x)}\,dx-\frac{b}{\gamma+1}\left(\int_\O \frac{|\nabla \Delta t\phi_0|^{p(x)}}{p(x)}\,dx\right)^{\gamma+1} - \lambda\int_\O \frac{f(x)}{q(x)}|t\phi_0|^{q(x)}\,dx \notag\\
&  & - \int_\O \frac{g(x)}{r(x)}|t\phi_0|^{r(x)}\,dx - \int_\O h(x)t\phi_0\,dx \notag\\
& \leq & at^{p_+}\int_\O\frac{|\nabla \Delta \phi_0|^{p(x)}}{p(x)}\,dx -\frac{bt^{p_-(\gamma+1)}}{\gamma+1}\left(\int_\O \frac{|\nabla \Delta \phi_0|^{p(x)}}{p(x)}\,dx\right)^{\gamma+1} + \frac{\lambda t^{q_+}}{q_-}\int_\O |f(x)||\phi_0|^{q(x)}\,dx \notag\\
& & - \frac{t^{r_-}}{r_-}\int_{\O_0} g(x)|\phi_0|^{r(x)}\,dx - t\int_\O h(x)\phi_0\,dx.
\end{eqnarray*}
Since $1< q_+ < p_+ < p_-(\gamma+1) < r_-$, we have $J(t\phi_0) \to -\infty$ as $t\to \infty$. So, for some $t_0>1$ large enough, we deduce that $\|t_0\phi_0\|_X > \rho$ and $J(t_0\phi_0) < 0$. Choosing $e = t_0\phi_0$, the proof of Lemma \ref{lem1} is completed. \qed

\begin{lem}\label{lem2}
Assume that $f$ satisfies $(H_1)-(H_2)$. Then the functional $J$ satisfies the Palais-Smale condition at level $c$ (popularly called $(PS)_c$ condition), where $c < \frac{\gamma a^\frac{\gamma+1}{\gamma}}{(\gamma+1)b^\frac{1}{\gamma}}$.
\end{lem}

{\bf Proof.} Let $\{u_n\}$ be a $(PS)_c$ sequence of $J$ such that $c < \frac{\gamma a^\frac{\gamma+1}{\gamma}}{(\gamma+1)b^\frac{1}{\gamma}}$, that is
\begin{eqnarray}\label{er1}
J(u_n) \to c, \quad J'(u_n) \to 0 \mbox{ in } X^*, \quad n \to \infty,
\end{eqnarray}
where $X^*$ is the dual space of $X$.

{\bf Step 1.} We first prove that $\{u_n\}$ is bounded in $X$. Arguing by contradiction, if $\{u_n\}$ is unbounded in $X$, up to a subsequence, we may assume that $\|u_n\|_X \to \infty$ as $n\to \infty$. Let $\theta$ be a fixed positive constant such that
$$
\theta \in \left(p_+, \min\left\{r_-, \frac{(p_-)^{\gamma+1} (\gamma+1)}{(p_+)^\gamma}\right\}\right).
$$

Then according to the conditions $(H_1)$ and $(H_2)$, for $n$ large enough, we have
\begin{eqnarray}\label{eer4}
c+1+\|u_n\|_X & \geq & J(u_n) - \frac{1}{\theta} \langle {J'(u_n), u_n} \rangle \notag\\
& = & a\int_\O \left(\frac{1}{p(x)}-\frac{1}{\theta}\right)|\nabla\Delta u_n|^{p(x)}\,dx - \frac{b}{\gamma+1}\left(\int_\O\frac{|\nabla\Delta u_n|^{p(x)}}{p(x)}\,dx\right)^{\gamma+1} \notag \\
& & + \frac{b}{\theta}\left(\int_\O\frac{|\nabla\Delta u_n|^{p(x)}}{p(x)}\,dx\right)^{\gamma}\int_\O |\nabla\Delta u_n|^{p(x)}\,dx + \lambda \int_\O \left(\frac{1}{\theta} - \frac{1}{q(x)}\right)f(x)|u_n|^{q(x)}\,dx \notag\\
& & + \int_\O \left(\frac{1}{\theta} - \frac{1}{r(x)}\right)g(x)|u_n|^{r(x)}\,dx + \left(\frac{1}{\theta}-1\right)\int_\O h(x)u_n\,dx \notag\\
& \geq & a\left(\frac{1}{p_+}-\frac{1}{\theta}\right)\|u_n\|^{p_+}_X + b\left(\frac{1}{\theta(p_+)^\gamma} - \frac{1}{(p_-)^{\gamma+1}(\gamma+1)}\right)\|u_n\|^{p_-(\gamma+1)}_X \notag\\
& & - \lambda C_1\left(\frac{1}{q_-}-\frac{1}{\theta}\right)|f|_{q_0(\cdot)}\|u_n\|^{q_+}_X - C_3\left(1-\frac{1}{\theta}\right)|h|_{\frac{p^*(\cdot)}{p^*(\cdot)-1}} \|u_n\|_X.
\end{eqnarray}
From \eqref{eer4}, it follows that
\begin{gather*}\label{er5}
c+1+\left[1+C_3\left(1-\frac{1}{\theta}\right)|h|_{\frac{p^*(\cdot)}{p^*(\cdot)-1}}\right] \|u_n\|_X + \lambda C_1\left(\frac{1}{q_-}-\frac{1}{\theta}\right)|f|_{q_0(\cdot)}\|u_n\|^{q_+}_X \geq a\left(\frac{1}{p_+}-\frac{1}{\theta}\right)\|u_n\|^{p_+}_X \notag\\
+ b\left(\frac{1}{\theta(p_+)^\gamma} - \frac{1}{(p_-)^{\gamma+1}(\gamma+1)}\right)\|u_n\|^{p_-(\gamma+1)}_X,
\end{gather*}
which is a contradiction since $\|u_n\|_X \to \infty$ as $n\to \infty$. So, $\{u_n\}$ is bounded in $X$ and the first assertion is proved.

{\bf Step 2.} Now, we prove that $\{u_n\}$ has a convergent subsequence in $X$. Indeed, by Proposition \ref{embedding}, the embedding $X \hookrightarrow L^{s(\cdot)}(\O)$ is compact, where $1\leq s(x)<p(x)^*$. Since $X$ is a reflexive Banach space, passing if necessary, to a subsequence, there exists $u\in X$ such that
\begin{equation}
\label{cvg}
u_n\rightharpoonup u\mbox{ in } X,\; u_n \to u \mbox{ in } L^{s(\cdot)}(\O),\;\; u_n(x)\to u(x), \mbox{ a.e. in } \O.
\end{equation}
From \eqref{deriv-fun}, we find that
\begin{eqnarray} \label{er'}
\langle J'(u_n), u_n-u\rangle &=& \left[a-b\left(\int_\O \frac{|\nabla \Delta u_n|^{p(x)}}{p(x)}\,dx\right)^{\gamma}\right]\int_\O|\nabla \Delta u_n|^{p(x)-2}\nabla\Delta u_n (\nabla\Delta u_n-\nabla\Delta u)dx \notag\\&\qquad-& \lambda\int_\O f(x)|u_n|^{q(x)-2}u_n(u_n-u)\,dx - \int_\O g(x)|u_n|^{r(x)-2}u_n(u_n-u)\,dx \notag\\&\qquad-& \int_\O h(x)(u_n-u)\,dx.
\end{eqnarray}
Meanwhile, by H\"older's inequality and \eqref{cvg} we estimate
\begin{eqnarray}\label{eqt0}
{  \left|\int_{\O}f(x)| u_n| ^{q(x)-2}u_n (u_n-u)\,dx\right|}   &\leq& \int_{\O}|f(x)| | u_n|^{q(x)-1}| u_n-u| \,dx\notag\\
&\leq& |f|_{q_0(\cdot)+\eta}{\Big|{|  u_n|} ^{q(x)-1}\Big|} _{\frac{p^*(\cdot)}{q(\cdot)-1}}| u_n-u| _{\theta(x)}\notag\\
&\leq& \max\left\{\|u\|^{q_+-1}_X, \|u\|^{q_--1}_X\right\}|f|_{q_0(\cdot)+\eta}| u_n-u| _{\theta_1(\cdot)},
\end{eqnarray}
where $\theta_1\in C(\overline{\O})$ such that $$\frac{1}{q_0(x)+\eta}+ \frac{q(x)-1}{p^*(x)}+\frac{1}{\theta_1(x)}=1.$$ We can easily verify that
$$\theta_1(x):=\frac{p^*(x)(p^*(x)+\eta p^*(x)-q(x)\eta)}{\eta{p^*(x)}^2+ (-2q(x)\eta+\eta+1)p^*(x)+q(x)\eta(q(x)-1)}< p^*(x).$$
So, thanks to \eqref{cvg} we can deduce that
\begin{equation}\label{eqt1}
| u_n-u| _{\theta_1(\cdot)}\to 0 \mbox{ as } n\to \infty.
\end{equation}
Combining this and the fact that  $\{u_n\}$ is bounded in $X$, we infer from \eqref{eqt0} and  \eqref{eqt1}  that
\begin{equation}\label{cvgu0}
\lim_{n\to \infty}\int_{\O}f(x)| u_n| ^{p(x)-2}u_n (u_n-u)\,dx=0.
\end{equation}
Similarly, we obtain
\begin{equation}\label{cvgu1}
\lim_{n\to \infty}\int_{\O}g(x)| u_n| ^{r(x)-2}u_n (u_n-u)\,dx=0 \mbox{ and } \lim_{n\to \infty}\int_{\O}h(x)(u_n-u)\,dx=0.
\end{equation}
By \eqref{er1}, we have
\[
\langle J'(u_n),u_n-u\rangle \to 0.
\]
So, from \eqref{cvgu0} and \eqref{cvgu1}, we can deduce that \eqref{er'} implies
\begin{equation}
\label{gj}
\left[a-b\left(\int_\O \frac{|\nabla \Delta u_n|^{p(x)}}{p(x)}\,dx\right)^{\gamma}\right]\int_\O|\nabla \Delta u_n|^{p(x)-2}\nabla\Delta u_n (\nabla\Delta u_n-\nabla\Delta u)\,dx\to 0.
\end{equation}
Since $\{u_n\}$ is bounded in $X$, passing to a subsequence, if necessary, we may assume that when $n\to \infty$
then
\[
\int_\O \frac{|\nabla \Delta u_n|^{p(x)}}{p(x)}\,dx\to t_0\geq 0.
\]
If $t_0 = 0$ then $\{u_n\}$ converges strongly to $u = 0$ in $X$ and the proof is finished. Otherwise, we need to consider the following two cases:
\\
{\bf Case 1.} If $t_0\neq \left(\frac ab\right)^\frac{1}{\gamma}$ then $a-b\left(\int_\O \frac{|\nabla \Delta u_n|^{p(x)}}{p(x)}dx\right)^{\gamma}\to 0$ is not true and no
subsequence of \\$\{a-b\left(\int_\O \frac{|\nabla \Delta u_n|^{p(x)}}{p(x)}\,dx\right)^{\gamma}\to 0\}$  converges to zero. Therefore, there exists $\delta > 0$
such that \\$\left| a-b\left(\int_\O \frac{|\nabla \Delta u_n|^{p(x)}}{p(x)}\,dx\right)^{\gamma}\right|  >\delta>0$ when $n$ is large enough. So, it is clear that
\begin{equation}
\left\{a-b\left(\int_\O \frac{|\nabla \Delta u_n|^{p(x)}}{p(x)}\,dx\right)^{\gamma}\to 0\right\}  \mbox{ is
bounded}.
\end{equation}
{\bf Case 2.} If $t_0 =\left(\frac ab\right)^\frac{1}{\gamma}$ then $a-b\left(\int_\O \frac{|\nabla \Delta u_n|^{p(x)}}{p(x)}dx\right)^{\gamma}\to 0$.
\\We define
\[
\varphi(u)=\lambda\int_\O \frac{f(x)}{q(x)}|u|^{q(x)}\,dx+\int_\O \frac{g(x)}{r(x)}|u|^{r(x)}\,dx+\int_\O h(x)u\,dx,\; \mbox{for all} \; u\in X.
\]
Then
\[
\langle \varphi'(u),v\rangle=\lambda\int_\O f(x)|u|^{q(x)-2}uv\,dx + \int_\O g(x)|u|^{r(x)-2}uv\,dx + \int_\O h(x)v\,dx,\; \mbox{for all} \; v\in X.
\]
It follows that
\[
\langle \varphi'(u_n)-\varphi'(u),v\rangle=\lambda\int_\O f(x)(|u_n|^{q(x)-2}u_n-|u|^{q(x)-2}u)v\,dx
+\int_\O g(x)(|u_n|^{r(x)-2}u_n-|u|^{r(x)-2}u)v\,dx.
\]
To complete the argument we need the following lemma.

\begin{lem}
\label{jkjk}
Let $u_n,u\in X$  be such that \eqref{cvg} holds. Then passing to a subsequence if necessary, the following properties hold:
\begin{enumerate}
  \item[(i)] $\lim_{n\to \infty}\int_\O f(x)(|u_n|^{q(x)-2}u_n-|u|^{q(x)-2}u)v\,dx=0$;
  \item[(ii)] $\lim_{n\to \infty}\int_\O g(x)(|u_n|^{r(x)-2}u_n-|u|^{r(x)-2}u)v\,dx=0$;
  \item[(iii)] $\langle \varphi'(u_n)-\varphi'(u),v\rangle\to 0,\;v\in X$.
\end{enumerate}
\end{lem}

{\bf Proof.} Let $\theta_2\in C(\overline{\O})$ be such that $\theta_2(x)<\frac{q(x)-1}{q(x)}p^*(x)$. From \eqref{cvg}, we deduce that $u_n\to u$ in $L^{\theta_2(\cdot)}(\O)$ which implies that
\begin{equation*}
\label{cvgi}
\left||u_n|^{q(x)-2}u_n- |u|^{q(x)-2}u\right|_{\frac{\theta_2(\cdot)}{q(\cdot)-1}}\to 0 \mbox{ as } n\to \infty.
\end{equation*}
Now, we define
$$q_0(x) = \frac{p^*(x)}{p^*(x)-q(x)},\qquad l_1(x)=\frac{\theta_2(x)}{q(x)-1},\qquad l_2(x)=\frac{p^*(x)\theta_2(x)}{q(x)\theta_2(x)-p^*(x)q(x)+p^*(x)},$$
where
$$l_2(x)<p^*(x), \quad\mbox{ and }\quad \frac{1}{q_0(x)}+\frac{1}{l_1(x)}+\frac{1}{l_2(x)}=1.$$
So, thanks to H\"older's inequality we can deduce that
\begin{eqnarray*}
\left|\int_\O\Big(f(x)|u_n|^{q(x)-2}u_n-f(x)|u|^{q(x)-2}u\Big)v\,dx\right|&\leq& \int_\O |f(x)|\Big||u_n|^{q(x)-2}u_n-|u|^{q(x)-2}u\Big||v|\,dx\nonumber\\
&\leq& |f|_{q_0(\cdot)}\Big||u_n|^{p(x)-2}u_n-|u|^{p(x)-2}u\Big|_{l_1(\cdot)}|v|_{l_2(\cdot)}\nonumber\\
&\leq& C_4 \Big||u_n|^{p(x)-2}u_n-|u|^{p(x)-2}u\Big|_{l_1(\cdot)}\|v\|_X\nonumber\\
&\to& 0.
\end{eqnarray*}
 By a slight modification of the proof above, we
can also prove assertion $(ii)$, so we omit the details.
\\
Finally, assertion $(iii)$ follows by combining parts $(i)$ and $(ii)$.
Consequently, $\|  \varphi'(u_n)-\varphi'(u)\| _{X^*}\to 0$ and $\varphi'(u_n)\to\varphi'(u)$. \qed

We can now complete the proof of Case 2. By Lemma \ref{jkjk} and since\\
 $\langle J'(u),v\rangle=\left[a-b\left(\int_\O \frac{|\nabla \Delta u|^{p(x)}}{p(x)}\,dx\right)^{\gamma}\right]\int_\O|\nabla \Delta u|^{p(x)-2}\nabla\Delta u \nabla \Delta v \,dx -\langle \varphi'(u),v\rangle$,\quad $\langle J'(u),v\rangle\to 0$\\
  and
   $$a-b\left(\int_\O \frac{|\nabla \Delta u|^{p(x)}}{p(x)}dx\right)^{\gamma}\to 0,$$  we can deduce that $\varphi'(u_n)\to 0\;(n\to\infty)$, i.e.,~\[
\langle \varphi'(u),v\rangle=\lambda\int_\O f(x)|u|^{q(x)-2}uv\,dx + \int_\O g(x)|u|^{r(x)-2}uv\,dx + \int_\O h(x)v\,dx,\; \mbox{for all} \; v\in X,
\]
 and therefore
\[
\lambda f(x)|u(x)|^{q(x)-2}u(x)+g(x)|u(x)|^{r(x)-2}u(x)+h(x)=0\mbox{ for a.e.} \; x \in \O.
\]
By the fundamental lemma of the variational method (see~\cite{Willem}) it
follows that $u=0$. Thus
\begin{gather*}
\varphi(u_n)=\lambda\int_\O \frac{f(x)}{q(x)}|u_n|^{q(x)}\,dx+\int_\O \frac{g(x)}{r(x)}|u_n|^{r(x)}\,dx+\int_\O h(x)u_n\,dx\\\to \lambda\int_\O \frac{f(x)}{q(x)}|u|^{q(x)}\,dx+\int_\O \frac{g(x)}{r(x)}|u|^{r(x)}\,dx+\int_\O h(x)u\,dx=0.
\end{gather*}
Hence, we can deduce that
\begin{eqnarray*}
J(u_n)&=&a\int_\O \frac{|\nabla \Delta u|^{p(x)}}{p(x)}\,dx-\frac{b}{\gamma+1}\left(\int_\O \frac{|\nabla \Delta u|^{p(x)}}{p(x)}\,dx\right)^{\gamma+1} - \lambda\int_\O \frac{f(x)}{q(x)}|u|^{q(x)}\,dx \notag\\
&  & - \int_\O \frac{g(x)}{r(x)}|u|^{r(x)}\,dx - \int_\O h(x)u\,dx\to \frac{\gamma a^\frac{\gamma+1}{\gamma}}{(\gamma+1)b^\frac{1}{\gamma}}.
\end{eqnarray*}
 This is a contradiction since $J(u_n)\to c<\frac{\gamma a^\frac{\gamma+1}{\gamma}}{(\gamma+1)b^\frac{1}{\gamma}}$, hence
  $$a-b\left(\int_\O \frac{|\nabla \Delta u_n|^{p(x)}}{p(x)}dx\right)^{\gamma}\to 0$$ is not true and similarly to Case $1$, we have that
\begin{equation*}
\left\{a-b\left(\int_\O \frac{|\nabla \Delta u_n|^{p(x)}}{p(x)}\,dx\right)^{\gamma}\to 0\right\}  \mbox{ is
bounded}.
\end{equation*}
So, it follows from  the two cases above that
\[
\int_\O|\nabla \Delta u_n|^{p(x)-2}\nabla\Delta u_n (\nabla\Delta u_n-\nabla\Delta u)\,dx\to 0.
\]
Applying $(S_+)$ mapping theory (see \cite{CP} for $r=3$), we can now deduce that $\| u_n\|_X \to\| u\|_X $ as $n\to \infty$, which means
that $J$ satisfies the $(PS)_c$ condition. This completes the proof.\qed

\begin{rem}
The $(PS)_c$ condition is not satisfied for $c>\frac{\gamma a^\frac{\gamma+1}{\gamma}}{(\gamma+1)b^\frac{1}{\gamma}}$.
\\Indeed,
\begin{eqnarray*}
J(u)&=&a\int_\O \frac{|\nabla \Delta u|^{p(x)}}{p(x)}\,dx-\frac{b}{\gamma+1}\left(\int_\O \frac{|\nabla \Delta u|^{p(x)}}{p(x)}\,dx\right)^{\gamma+1} - \lambda\int_\O \frac{f(x)}{q(x)}|u|^{q(x)}\,dx\\
&\quad-& \int_\O \frac{g(x)}{r(x)}|u|^{r(x)}\,dx - \int_\O h(x)u\,dx\\
&\leq&a\int_\O \frac{|\nabla \Delta u|^{p(x)}}{p(x)}\,dx-\frac{b}{\gamma+1}\left(\int_\O \frac{|\nabla \Delta u|^{p(x)}}{p(x)}\,dx\right)^{\gamma+1}\leq \frac{\gamma a^\frac{\gamma+1}{\gamma}}{(\gamma+1)b^\frac{1}{\gamma}},
\end{eqnarray*}
So, if $\{u_n\}$ is a $(PS)_c$ sequence of $J$, then we have $c \leq \frac{\gamma a^\frac{\gamma+1}{\gamma}}{(\gamma+1)b^\frac{1}{\gamma}}$, which is a contradiction.
\end{rem}

\section{\bf  Proof of Theorem \ref{the1.1}}\label{sect4}
In view of Lemmas \ref{lem1}, \ref{lem2} and the Mountain pass theorem in \cite{Willem}, there exists a weak solution $u_1$ of problem \eqref{high} with $\overline c = J(u_1)>0$. We will show that there exists a second weak solution $u_2 \ne u_1$ by using the Ekeland variational principle. First, let us choose $\psi_0 \in C_0^\infty(\O)$ such that $\int_\O h(x)\psi_0dx >0$. Now, for all $t \in (0,1)$ small enough, we have
\begin{eqnarray}\label{er4}
J(t\psi_0) & = & a\int_\O \frac{|\nabla \Delta t\psi_0|^{p(x)}}{p(x)}\,dx-\frac{b}{\gamma+1}\left(\int_\O \frac{|\nabla \Delta t\psi_0|^{p(x)}}{p(x)}\,dx\right)^{\gamma+1} - \lambda\int_\O \frac{f(x)}{q(x)}|t\psi_0|^{q(x)}\,dx \notag\\
&  & - \int_\O \frac{g(x)}{r(x)}|t\psi_0|^{r(x)}\,dx - \int_\O h(x)t\psi_0\,dx \notag\\
& \leq & at^{p_-}\int_\O\frac{|\nabla \Delta \phi_0|^{p(x)}}{p(x)}\,dx -\frac{bt^{p_+(\gamma+1)}}{\gamma+1}\left(\int_\O \frac{|\nabla \Delta \phi_0|^{p(x)}}{p(x)}\,dx\right)^{\gamma+1} + \frac{\lambda t^{q_-}}{q_-}\int_\O |f(x)||\phi_0|^{q(x)}\,dx \notag\\
& & + \frac{t^{r_-}}{r_-}\int_{\O_0} g(x)|\phi_0|^{r(x)}\,dx - t\int_\O h(x)\phi_0\,dx < 0
\end{eqnarray}
since
$$
at^{p_-}\int_\O\frac{|\nabla \Delta \phi_0|^{p(x)}}{p(x)}\,dx + \frac{\lambda t^{q_-}}{q_-}\int_\O |f(x)||\phi_0|^{q(x)}\,dx - t\int_\O h(x)\phi_0\,dx < 0, \quad 1 < q_- < p_-
$$
and
$$
\frac{t^{r_-}}{r_-}\int_{\O_0} g(x)|\phi_0|^{r(x)}\,dx - \frac{bt^{p_+(\gamma+1)}}{\gamma+1}\left(\int_\O \frac{|\nabla \Delta \phi_0|^{p(x)}}{p(x)}\,dx\right)^{\gamma+1} < 0, \quad p_+(\gamma+1) < r_-
$$
for all $t\in (0,1)$ small enough.

By Lemma \ref{lem1}, it follows that on the boundary of the ball centered at the origin and of radius $\rho$ in $X$, denoted by $B_\rho(0)$, we have
$$
\inf_{u\in\partial B_\rho(0)}J(u)>0.
$$

On the other hand, again by Lemma \ref{lem1}, the functional $J$ is bounded from below on $B_\rho(0)$. Moreover, by \eqref{er1}, there exists $\psi_0 \in X$ such that $J(t\psi_0)<0$ for all $t>0$ small enough. It follows that
$$
-\infty < \underline c=\inf_{u\in \overline B_\rho(0)}J(u)<0.
$$

Let us choose $\varepsilon > 0$ such that $0 < \varepsilon<\inf_{u \in \partial B_\rho(0)} J(u)- \inf_{u \in \overline B_\rho(0)}J(u)$. Applying Ekeland's variational principle \cite{Ekeland} to the functional $J: \overline B_\rho(0)  \to \R$, it follows that there exists $u_\varepsilon \in \overline B_\rho(0)$ such that
\begin{align*}
J(u_\varepsilon) & < \inf_{u \in \overline B_\rho(0)}J(u)+\varepsilon, \\
J(u_\varepsilon) & < J(u) + \varepsilon\|u-u_\varepsilon\|_X, \quad u \ne u_\varepsilon,
\end{align*}
so we have $J(u_\varepsilon) <  \inf_{u \in \partial B_\rho(0)}J(u)$ and thus, $u_\varepsilon \in B_\rho(0)$.

Now, we define the functional $I: \overline B_\rho(0) \to \R$ by $I(u) = J(u) + \varepsilon \|u-u_\varepsilon\|_X$. It is clear that $u_\varepsilon$ is a minimum point of $I$ and thus
$$
\frac{I(u_\varepsilon + \tau v)-I(u_\varepsilon)}{t} \geq 0
$$
for all $\tau > 0$ small enough and all $v \in B_\rho(0)$. The above information shows that
$$
\frac{J(u_\varepsilon+\tau v)-J(u_\varepsilon)}{\tau }+\varepsilon \|v\|_X \geq 0.
$$
Letting $\tau \to 0^+$, we deduce that
$$
\left\langle {J'(u_\varepsilon),v} \right\rangle + \varepsilon \|v\|_X \geq 0,
$$
which leads to $\|J'(u_\varepsilon)\|_{X^\ast} \leq \varepsilon$. Therefore, there exists a sequence $\{u_n\} \subset B_\rho(0)$ such that
\begin{equation}\label{e3.27}
J(u_n) \to \underline c  = \inf_{u \in \overline B_\rho(0)}J(u) < 0 \text{ and } J'(u_n) \to 0 \text{ in } X^\ast \text{ as } n \to \infty.
\end{equation}

By Lemma \ref{lem2}, the sequence $\{u_n\}$ converges strongly to some $u_2$ as $n \to \infty$. Moreover, since $J \in C^1(X,\R)$, by \eqref{e3.27} it follows that $J'(u_2) = 0$. Thus, $u_2$ is a nontrivial weak solution of problem \eqref{high} with negative energy $J(u_2) = \underline c<0$.

Finally, we point out the fact that $u_1 \ne u_2$ since $J(u_1) = \overline c > 0 > \underline c = J(u_2)$. The proof of Theorem \ref{the1.1} is now complete.\qed

\begin{rem}
The proof of Theorem \ref{the1.2} is very similar.
\end{rem}

\subsection*{\bf Acknowledgments}
The first author was supported by the Tunisian Military Research Center
for Science and Technology Laboratory LR19DN01. The second author was supported by Vietnam National Foundation for Science and Technology Development (NAFOSTED) (Grant N.101.02.2017.04). The third author was supported by the Slovenian Research Agency grants P1-0292, N1-0114, N1-0083, N1-0064, and J1-8131. We gratefully acknowledge the referees for their comments and suggestions.




\begin{thebibliography}{777}
\bibliographystyle{alpha}
\bibitem{AM1}
E. Acerbi, G. Mingione, \emph{Gradient estimates for the $p(x)$-Laplacean system,} J. Reine Angew. Math., \textbf{584}, (2005) 117-148.

\bibitem{AMC14} G.A. Afrouzi, M. Mirzapour, N.T. Chung, \emph{Existence and multiplicity of solutions for Kirchhoff type problems involving $p(x)$-biharmonic operators}, Z. Anal. Anwend., $\mathbf{33}$, (2014) 289-303.



\bibitem{BR16}
S. Baraket, V.D. R\u{a}dulescu, \emph{Combined effects of concave-convex nonlinearities in a fourth-order problem with variable exponent}, Adv. Nonlinear Stud., $\mathbf{16}$(3), (2016): 409.

\bibitem{BRR16}
M.M. Boureanu, V.D. R\u{a}dulescu, D.D. Repov\v{s}, \emph{On a $p(\cdot)$-biharmonic problem with no-flux boundary condition}, Comput. Math. Appl., $\mathbf{72}$ (9), (2016) 2505-2515.

\bibitem{JF}
 J. Chabrowski, Y. Fu, \emph{Existence of solutions for $p(x)$-Laplacian problems on a bounded domain}, J. Math. Anal. Appl., $\mathbf{306}$, (2005) 604-618.

\bibitem{CHL}
C. Chen, J. Huang, L. Liu, .\emph{Multiple solutions to the nonhomogeneous p-Kirchhoff elliptic equation with concave-convex nonlinearities}, Appl. Math. Lett., $\mathbf{26}$(7), (2013) 754-759.

\bibitem{CKW}
C.Y. Chen, Y.C. Kuo, T.F. Wu,  \emph{The Nehari manifold for a Kirchhoff type problem involving sign-changing weight functions}, J. Differential Equations $\mathbf{250}$, (2011) 1876–1908.

\bibitem{CLR}
 Y. Chen, S. Levine, M. Rao, \emph{Variable exponent, linear growth functionals in image restoration}, SIAM J. Appl. Math., $\mathbf{66}$, (2006) 1383-1406.

\bibitem{Chung2012}
N.T. Chung, \emph{Multiplicity results for a class of $p(x)$-Kirchhoff type equations with combined nonlinearities}, Elec. J. Qual. Theory Diff. Equ., $\mathbf{2012}$(42), (2012) 1-13.

\bibitem{Chung13}
N.T. Chung, \emph{Multiple solutions for a $p(x)$-Kirchhoff-type equation with sign-changing nonlinearities}, Complex Var. Elliptic Equ., $\mathbf{58}$(12) (2013), 1637-1646.

\bibitem{C182}
N.T. Chung, \emph{Some remarks on a class of $p(x)$-Laplacian Robin eigenvalue problems}, Mediterr. J. Math., $\mathbf{15}$(4), (2018): 147.

\bibitem{C181}
 N.T. Chung, \emph{Existence of solutions for perturbed fourth order elliptic equations with variable exponents}, Electron. J. Qual. Theory Differ. Equ., $\mathbf{2018}$(96), (2018) 1-19.
 
\bibitem{CN10}
N.T. Chung, Q.A. Ngo, \emph{Multiple solutions for a class of quasilinear elliptic equations of $p(x)$-Laplacian type with nonlinear boundary conditions}, Proc. Royal Soc. Edinburgh Sect. A: Mathematics, $\mathbf{140}$(2), (2010) 259-272.

\bibitem{CP}
 F. Colasuonno, P. Pucci, \emph{Multiplicity of solutions for $p(x)$-polyharmonic Kirchhoff equations}, Nonlinear Anal., $\mathbf{74}$, (2011) 5962-5974.

\bibitem{DH}
G. Dai and R. Hao, \emph{Existence of solutions for a $p(x)-$Kirchhoff-type equation}, J. Math. Anal. Appl., $\mathbf{359}$, (2009) 275-284.

\bibitem{Diening}
L. Diening, P. Harjulehto, P. H\"ast\"o, M. R\r{u}\v{z}i\v{c}ka, \emph{Legesgue and Sobolev spaces with  variable exponents}, Lecture Notes in Mathematics  $\mathbf{2017}$, Springer-Verlag, Heidelberg, 2011.

\bibitem{Ekeland}
I. Ekeland, \emph{On the variational principle}, J. Math. Anal. Appl., $\mathbf{47}$, (1974), 324-353.

\bibitem{AMM}
A. El Amrouss, F. Moradi, M. Moussaoui, \emph{Existence of solutions for fourth-order PDEs with variable exponents}, Electron. J. Differ. Equ., $\mathbf{2009}$, (153), (2009) 1-13.

\bibitem{FS}
X. L. Fan, J. S. Shen, D. Zhao, \emph{Sobolev embedding theorems for spaces $W^{k,p(x)}(\O)$}, J. Math. Anal. Appl., $\mathbf{262}$, (2001) 749-760.

\bibitem{FZ1}
 X. L. Fan, D. Zhao, \emph{On the spaces $L^{p(x)}$ and $W^{m,p(x)}$}, J. Math. Anal. Appl., $\mathbf{263}$, (2001) 424-446.

\bibitem{H}
T. C. Halsey, \emph{Electrorheological fluids,} Science, $\mathbf{258}$, (1992) 761-766.

\bibitem{H2}
M.K. Hamdani, \emph{On a nonlocal asymmetric Kirchhoff problems,} Asian-European J. Math., (2019), doi: 10.1142/S1793557120300018

\bibitem{HCA}
M.K. Hamdani, N.T. Chung, M.B. Aminlouee, \emph{Infinitely many solutions for a new class of Schr\"{o}dinger-Kirchhoff type equations in $\R^N$ involving the fractional $p$-Laplacian,} J. Elliptic Parabol. Equ., doi: 10.1007/s41808-020-00093-7.

\bibitem{HHMR}
M.K. Hamdani, A. Harrabi, F. Mtiri, and D.D. Repov\v{s}, \emph{Existence and multiplicity results for a new $p(x)-$Kirchhoff problem}.  Nonlinear Anal., $\mathbf{190}$ (2020): 111598.

\bibitem{HR}
M. K. Hamdani, D.D. Repov\v{s}, \emph{Existence of solutions for systems arising in electromagnetism}, J. Math. Anal. Appl., $\mathbf{486}$(2) (2020):123898.

\bibitem{HZCR}
M.K. Hamdani, J. Zuo, N. T. Chung, D.D. Repov\v{s}, \emph{Multiplicity of solutions for a class of fractional $ p(x,\cdot)$-Kirchhoff-type problems without the Ambrosetti–Rabinowitz condition,} Bound. Value Probl., $\mathbf{2020}$, (2020):150.

\bibitem{HCX}
J. C. Huang, C. S. Chen, Z. H. Xiu, \emph{Existence and multiplicity results for a $p$-Kirchhoff equation with a concave-convex term}, Appl. Math. Lett., $\mathbf{26}$ (2013) 1070-1075.

\bibitem{KR}
O. Kov\'a\v{c}ik, J. R\'akosn\'ik, \emph{On spaces $L^{p(x)}$ and $W^{k,\,p(x)}$}, Czechoslovak Math. J., $\mathbf{41}$, (1991) 592-618.

\bibitem{LLP}
F. Li, Z. Li, L. Pi, \emph{Variable exponent functionals in image restoration}, Appl. Math. Comput., $\mathbf{216}$ (3), (2010) 870-882.

\bibitem{LMZ}
Y.X. Li, M. Mei, K.J. Zhang, \emph{Existence of multiple nontrivial solutions for a $p$-Kirchhoff type elliptic problem involving sign-changing weight functions}, Discrete Contin. Dyn. Syst., Ser. B, $\mathbf{21}$, (2016) 883-908.

\bibitem{LRDZ}
G. Li, V. D. R\u{a}dulescu, D. D. Repov\v{s}, Q. Zhang, \emph{Nonhomogeneous Dirichlet problems without the Ambrosetti-Rabinowitz condition}, Topol. Methods Nonlinear Anal. $\mathbf{51}$ (1), (2018) 55-77.

\bibitem{MAE}
R.A. Mashiyev, H. Alisoy, I. Ekincioglu, \emph{Existence of one weak solution for $p(x)$-biharmonic equations involving a concave-convex nonlinearity}, Matemati\v cki Vesnik, $\mathbf{69}$, (2017) 296-307.

\bibitem{MTTL}
M. Massar, M. Talbi, N. Tsouli, H. Lebrimchi, \emph{On $p(x)$-Kirchhoff equations with concave-convex terms in unbounded domains}, J. Nonlinear Funct. Anal. $\mathbf{2018}$ (2018), https://doi.org/10.23952/jnfa.2018.10

\bibitem{MR06}
 M. Mih\v{a}ilescu, V. D.  R\u{a}dulescu, \emph{A multiplicity result for a nonlinear degenerate problem arising in the theory of electrorheological fluids}, Proc. R. Soc. A $\mathbf{462}$, (2006) 2625-2641.

\bibitem{MR10}
 M. Mih\v{a}ilescu, V. D.  R\u{a}dulescu, \emph{Eigenvalue problems with weight and variable exponent for the Laplace operator}, Anal. Appl., $\mathbf{8}$, (2010) 235-246.

\bibitem{R15}
V. D. R\u{a}dulescu, \emph{Nonlinear elliptic equations with variable exponent: old and new}, Nonlinear Anal., $\mathbf{121}$ (2015) 336-369.

\bibitem{RD}
V. D. R\u{a}dulescu, D. D. Repov\v{s}, \emph{Partial differential equations with variable exponents: variational methods and qualitative analysis}, CRC Press, Boca Raton, 2015.

\bibitem{Rahal}
 B. Rahal, \emph{Existence results of infinitely many solutions for $p(x)-$Kirchhoff type triharmonic operator with Navier boundary conditions,} J. Math. Anal. Appl., $\mathbf{478}$, (2019) 1133-1146.

\bibitem{RR2}
 K.R. Rajagopal, M. R\r{u}\u{z}i\u{c}ka; \emph{Mathematical modeling of electrorheological materials,} Contin. Mech. Thermodyn., $\mathbf{13}$, (2001) 59-78.

\bibitem{R3}
M. R\r{u}\u{z}i\u{c}ka; \emph{Electro-rheological fluids: modeling and mathematical theory,} Lecture Notes in Math.
$\mathbf{1784}$, Springer, Berlin, (2000).

\bibitem{Willem}
M. Willem, \newblock {\em Minimax theorems}, {\em Progress in Nonlinear Differential Equations and their Applications} $\mathbf{24}$, \newblock Birkh\"{a}user Boston, Inc., Boston, MA, 1996.

\bibitem{Y}
Z. Yucedag, \emph{Existence of solutions for $p(x)$ Laplacian equations without Ambrosetti-Rabinowitz type condition}, Bull. Malay. Math. Sci. Soc., $\mathbf{38}$ (3), (2015) 1023-1033.

\bibitem{Z}
V. V. Zhikov, \emph{Averaging of functionals of the calculus of variations and elasticity theory,} Math. USSR. Izv, $\mathbf{29}$, (1987) 33-66.

\bibitem{Zhou}
Z. Zhou, \emph{On a $p(x)$-biharmonic problem with Navier boundary condition}, Bound. Value Probl., $\mathbf{2018}$, (2018): 149.

 \end{thebibliography}
\end{document}